\documentclass[12pt,a4]{article}
\usepackage{amssymb,amsmath}

\title{Dissipative hydrodynamic models for the diffusion of impurities in a gas\thanks{Support by the European network HYKE, funded by the EC as
contract HPRN-CT-2002-00282, is acknowledged.}}

\author{Stephane Brull\thanks{LATP. Centre de math\'ematiques et d'informatique. Universit\'e de Provence, 39 rue Joliot-Curie, 13453 Marseille. France. Email: {\tt brull@cmi.univ-mrs.fr}}\and Lorenzo Pareschi\thanks{Department of Mathematics, University of
Ferrara, Via Machiavelli 35, 44100 Ferrara, Italy. E-mail: {\tt
pareschi@dm.unife.it }}}

\begin{document}

\maketitle

\begin{abstract}
Recently linear dissipative models of the Boltzmann equation have
been introduced in \cite{ST, LT}. In this work, we consider the
problem of constructing suitable hydrodynamic approximations for
such models.
\end{abstract}

\section{Introduction}
The dissipative linear Boltzmann equation describes the dynamic of
a set of particles with mass $m_1$ interacting inelastically with
a background gas in thermodynamical equilibrium composed of
particles with mass $m \ll m_1$. For example, the case of fine
polluting impurities interacting with air or another gas is investigated in \cite{GM}.

As observed in \cite{ST}, the only conserved quantity is the number
of inelastic particles and as a result, a conventional hydrodynamic
approach of Euler type leads to a single equation describing the
advection (or advection-diffusion at the Navier-Stokes order) of
inelastic particles at the velocity of the background.

The aim of this note is to find hydrodynamic models for such
Boltzmann equation which posses equations for the momentum and the
temperature of the gas. Here, we present a closed set of
dissipative Euler equations for a pseudo-Maxwellian case which
generalizes the one considered in \cite{ST}.

Let us mention that the problem of finding suitable hydrodynamics
for inelastic interacting gases has been studied recently by
several authors (see \cite{BCG, BCGP, BST, PT} and the references
therein).


The paper is organized as follows. Section 2 deals with the linear
dissipative Boltzmann model and the pseudo-Maxwellian
approximation. Section 3 is devoted to discuss the problem of the
closure for the moment equations and to derive a dissipative Euler
system.

\section{The dissipative linear Boltzmann equation}
We consider the dissipative linear Boltzmann equation
\begin{equation} \label{Bo}
\frac{\partial f}{\partial t}(t,x,v) + v \cdot \nabla_{x}f(t,x,v) = Q(f)(t,x,v),\end{equation}
with,
\begin{equation} \label{Q}
 \frac{1}{2 \pi \lambda} \int_{ \mathbb{R}_{v}^{3} \times S^{2}} B(v,w,n) [ \frac{1}{e^{2}} f(v_{\ast}) M_{1}(w_{\ast}) - f(v) M_{1}(w)] dw dn.
\end{equation}
Here, $B(v,w,n)$ denotes the collision kernel, $\lambda$ the mean
free path and $e$ the restitution coefficient with $0 <e <1$. The
case $e=1$ corresponds to the elastic collision mechanism.

For the hard spheres model, the particles are assumed to be ideally
elastic balls and the corresponding collision kernel is given by
\begin{equation} \label{kernel}
 B(v,w,n) = |q | ,
\end{equation}
with $q = v-w$. The background is assumed to be in thermodynamic
equilibrium with given mass velocity $u_{1}$ and temperature
$T_{1}$ i.e. its distribution function $M_{1}$ is the normalized
Maxwellian given by
\begin{equation} \label{max}
M_{1}(v) = \frac{\rho_{1}}{(2 \pi T_{1})^{\frac{3}{2}}} \exp \big(
- \frac{(v - u_{1})^{2}}{2 T_{1}} \big).
\end{equation}
 Mass ratio and inelasticity are described by the following dimensionless parameters,
\begin{equation} \label{m}
\alpha = \frac{m_{1}}{m_{1} + m} \qquad \mbox{and} \qquad \beta =
\frac{1-e}{2},
\end{equation}
where $0 < \alpha < 1$ and $0 < \beta < \frac{1}{2}$.

In these conditions, it's possible to prove (see \cite{ST, LT})
that the stationary equilibrium states of the collision operators
are given by the Maxwellian distributions
\begin{equation}
    M^{\sharp}({v})=\left(\frac{m}{2\,\pi\, T^{\sharp}}\right)^{3/2}\exp\left\{-\frac{m({v}-u_{1})^{2}}{2\, T^{\sharp}}\right\}, \qquad {v}\in \mathbb{R}^{3}\,,
\label{eq-max}
\end{equation}
having the same mean velocity of the background and temperature
\begin{equation}
    T^{\sharp}=\frac{\left(1-\alpha\right) \left(1-\beta\right)}{1-\alpha\left(1-\beta\right)}T_{1}
\end{equation}
lower than the background one.

Here by analogy with \cite{BCG}, we consider an approximation of
the hard sphere model characterized by the assumption
\begin{equation} \label{ap}
|v - w| \simeq S(t,x),
\end{equation}
where $S(t,x)$ is a suitable function which takes into account the
fact that we have large relaxation rates for $|v-w|$ large and
small relaxation rates for $|v-w|$ small. Clearly since $v$ is
distributed accordingly to $f$ and $w$ accordingly to $M_1$ the
function $S$ cannot be simply a function of the temperature of a
single gas as in \cite{BCG}.

On the other hand, $M_1$ is given by (\ref{max}) and thus a
possible choice here consists in taking the expected value for
$|v-w|$ as choice of $S$. This gives
\begin{equation}
S(x,t)=\int_{ \mathbb{R}_{v}^{3}}\int_{ \mathbb{R}_{v}^{3}}
|v-w|f(x,v,t)M_1(x,w)\,dw\,dv=\int_{ \mathbb{R}_{v}^{3}} Z(x,v)
f(v) dv,
\end{equation} with
\begin{equation}
Z(x,v)=\int_{ \mathbb{R}_{v}^{3}} |v-w|M_1(x,w)\,dw.
\end{equation}
Off course simpler choices can be done. For example, similarly to
the case of a single gas, taking $S(x,t)=\mu\sqrt{T_r(x,t)}$ for a
suitable constant $\mu$, where $T_r$ is the normalized relative
``temperature'' given by
\[
T_r(x,t)=\frac{1}{3\rho(x)}\int f(x,v,t)|v-u_1(x)|^2\,dv.
\]
Note that at variance with \cite{BCG} here the ``temperature''
$T_r$ of the inelastic gas is measured with respect to the mean
velocity of the background. Thus only asymptotically for large
times it will correspond to the physical temperature.

Therefore this pseudo-Maxwellian model is given by
\begin{equation} \label{Bod}
\frac{\partial f}{\partial t} + v \cdot \nabla_{x}f =
\frac{S(t,x)}{2 \pi \lambda} \int_{ \mathbb{R}_{v}^{3} \times
S^{2}} [ \frac{1}{e^{2}} f(v_{\ast}) M_{1}(w_{\ast}) - f(v)
M_{1}(w)] dw dn.
\end{equation}

The above model represents a better approximation of the hard
sphere model with respect to the Maxwellian model considered in
\cite{ST} which corresponds simply to $S(x,t)= const$.


\section{Hydrodynamic limit and the Euler equation.}
To avoid the term $\frac{1}{e^{2}}$ in (\ref{Bod}), it is useful to
consider the weak form of (\ref{Bod}). More precisely, let us define
with $< \cdot , \cdot >$ the inner product in
$L^{1}(\mathbb{R}^{3})$. Given any regular test-function $\varphi(v)$,
it holds that
\begin{equation} \label{mo}
< \varphi, Q(f)> =  \frac{S(t,x) }{\lambda \pi}\int_{
\mathbb{R}_{w}^{3}} \int_{ \mathbb{R}_{v}^{3} \times S^{2}}
(\varphi(v^{\ast}) - \varphi(v)) f(v) M_{1}(w) dw dn,
\end{equation}
where the post-collisional velocity $v^{\ast}$ is defined by

\begin{equation} \label{col}
v^{\ast} = v -2 \alpha (1- \beta )( q \cdot n) n.
\end{equation}
Clearly $\varphi = 1$ is a collision invariant whereas $\varphi = v$ and $\varphi = v^{2}$ are not.

The existence of a Maxwellian equilibrium at non-zero temperature
(\ref{eq-max}) allows to construct hydrodynamic models for the
considered granular flow. However, here only the mass of the
inelastic particles is preserved. Thus the mass $\rho$ is the
unique hydrodynamic variable and the Euler system is reduced to
the single advection equation \cite{ST}
\begin{equation}
    \frac{\partial \rho}{\partial t}+\nabla\cdot(\rho
    u_1)=0.
\end{equation}

At this point, in order to perform a closure for the moment
equations such that the equations for the mean velocity and the
temperature of particles are preserved we assume the distribution
function $f$ to be the local Maxwellian at the mean velocity and
temperature of the gas
\begin{equation} \label{M}
M(x,v,t) = \frac{\rho (x,t)}{(2 \pi T(x,t))^{\frac{3}{2}}} \exp
\big( - \frac{(v - u(x,v))^{2}}{2 T(x,t)} \big).
\end{equation}

Taking $\varphi = v$ in (\ref{mo}) leads to
\begin{equation} \label{m2}
 < v , Q(M)> =  \frac{-2 \alpha (1-\beta) S(t,x)}{\lambda \pi} \int_{ \mathbb{R}_{w}^{3}} \int_{ \mathbb{R}_{v}^{3}}M(v) M_{1}(w) (\int_{ S^{2}} (q \cdot n)  n dn) dw dv .
\end{equation}
Following (\cite{BST},\cite{LT}), we get
\begin{equation} \label{dt}
\int_{ S^{2}} (q \cdot n)  n dn = \frac{4\pi}{3}  q .
\end{equation}
So, (\ref{m2}) has the following expression

\begin{equation} \label{v}
  < v , Q(M)> =  \frac{-8 \pi \alpha (1-\beta)S(t,x) }{3\lambda}  \int_{ \mathbb{R}_{w}^{3}} \int_{ \mathbb{R}_{v}^{3}}M(v) M_{1}(w) (v-w) dw dv .
\end{equation}
As,
\begin{equation}
\int_{ \mathbb{R}_{v}^{3}}vM(v)dv = \rho u \qquad \mbox{and} \qquad \int_{ \mathbb{R}_{v}^{3}}w M_{1}(w)dv = \rho_{1} u_{1},
\end{equation}
the first moment equation has the expression,

\begin{equation} \label{vi}
\frac{\partial}{\partial t}u + \nabla\cdot(\rho u \otimes u) +
\nabla_{x}(\rho T) = \frac{-4 S(t,x) \alpha (1-\beta)}{3\lambda}
\rho \rho_{1} (u -u_{1}).
\end{equation}

\noindent For the second moment, let us compute (\ref{mo}) with $\varphi = \frac{1}{2} |v|^{2}$. Hence,
\begin{eqnarray} \label{to}
< \frac{1}{2}|v|^{2} , Q(M)> \nonumber &=&  \frac{S(t,x) }{\lambda
\pi}\int_{ \mathbb{R}_{w}^{3}} \int_{ \mathbb{R}_{v}^{3} \times
S^{2}} \big[ -2\alpha (1-\beta) (q \cdot n) (v \cdot n)
\\& +&4 \alpha^{2} (1-\beta)^{2} |q \cdot n|^{2}  \big] M(v) M_{1}(w) dw dv dn. \qquad
\end{eqnarray}

\noindent Reasoning as in (\cite{BST},\cite{LT}), it holds that

\begin{eqnarray} \label{rt}
\int_{ S^{2}} |q \cdot n|^{2} dn = \frac{2\pi}{3} |q|^{2},
\\ \int_{ S^{2}} (q \cdot n)(v \cdot n) dn = \frac{2\pi}{3}(q \cdot n).
\end{eqnarray}

\noindent So, integrating the right-hand side of (\ref{to}) with
respect to the $n$ variable and using (\ref{rt}) leads to

\begin{eqnarray} \label{ti}
\nonumber  < \frac{1}{2}|v|^{2} , Q(M)> &=& \frac{-4  \alpha
(1-\beta)S(t,x)}{3\lambda} \int_{ \mathbb{R}_{w}^{3}} \int_{
\mathbb{R}_{v}^{3}}   (q \cdot v)  M(v) M_{1}(w) dw dv \qquad
\qquad
\\ &+& \frac{4  \alpha^{2} (1-\beta)^{2} S(t,x)}{3\lambda} \int_{ \mathbb{R}_{w}^{3}} \int_{ \mathbb{R}_{v}^{3}}  |q|^{2}   M(v) M_{1}(w) dw dv.\qquad \qquad
\end{eqnarray}

\noindent  As, $q = v-w$,

\begin{equation} \label{p}
 \int_{ \mathbb{R}_{w}^{3}} \int_{ \mathbb{R}_{v}^{3}}  |q|^{2} M(v) M_{1}(w) dw dv =  \int_{ \mathbb{R}_{w}^{3}} \int_{ \mathbb{R}_{v}^{3}}  (|v|^{2} - 2 v\cdot w + |w|^{2})   M(v) M_{1}(w) dw dv.
\end{equation}

\noindent As,

\begin{eqnarray} \label{pi}
\nonumber \frac{1}{2}  \int_{ \mathbb{R}_{v}^{3}} | v|^{2} M(v) dv = \rho (\frac{1}{2} |u|^{2} + \frac{3}{2} T),
\end{eqnarray}

\noindent it follows that

\begin{equation} \label{i1}
 \int_{ \mathbb{R}_{w}^{3}} \int_{ \mathbb{R}_{v}^{3}}  |q|^{2} M(v) M_{1}(w) dw dv =  \rho \rho_{1} (3 T + 3T_{1}+ |u|^{2} + |u_{1}|^{2}- 2 u_{1}\cdot u).
\end{equation}

\noindent  and

\begin{equation} \label{i2}
 \int_{ \mathbb{R}_{w}^{3}} \int_{ \mathbb{R}_{v}^{3}}  (|v|^{2} - v
 \cdot w)  M(v) M_{1}(w) dw dv =   \rho \rho_{1} (3 T + |u|^{2} - u \cdot u_{1}).
\end{equation}
By (\ref{i1}) and (\ref{i2}), the right-hand side of (\ref{ti}) is equal to

\begin{equation} \label{i}
2  \alpha^{2} (1-\beta)^{2}  \rho \rho_{1} (3 T + 3T_{1}+ |u|^{2} + |u_{1}|^{2} - 2 u_{1}\cdot u) -2\alpha (1-\beta) \rho \rho_{1} (3 T + |u|^{2} - u \cdot u_{1}).
\end{equation}
Finally, the left-hand side of (\ref{ti}) being computed by (\ref{pi}), we find the following dissipative Euler system
\begin{eqnarray}
\nonumber \frac{\partial \rho}{\partial t} + \nabla\cdot(\rho u)
&=& 0,
\\  \frac{\partial}{\partial t}u + \nabla\cdot (\rho u \otimes u) +
\nabla_{x}(\rho T) &=& \frac{4 \pi S(t,x) \alpha (1-\beta) }{3
\lambda} \rho \rho_{1} (u_1 -u)
\\ \nonumber \frac{\partial}{\partial t} \big( \rho u (
\frac{1}{2}|u|^{2} +\frac{3}{2}T ) \big) + \nabla\cdot \big( \rho
( \frac{1}{2}|u|^{2} +\frac{5}{2}T ) \big)&=&\frac{4 \pi \rho
\rho_{1} S(t,x)}{3\lambda} D(x,t)
\end{eqnarray}
where
\begin{equation}
D(x,t)=\alpha^{2} (1-\beta)^{2}
 (3 T + 3T_{1}+ |u|^{2} + |u_{1}|^{2} - 2 u_{1}\cdot u) - \alpha (1-\beta)   (3 T + |u|^{2} - u \cdot u_{1}).
\end{equation}

\section{Conclusion}
We derived hydrodynamic approximations for linear dissipative
Boltzmann equations that keep the equations for the mean velocity
and the temperature of particles. To this aim the closure of the
moment system is performed with respect to a local Maxwellian
state which is not an equilibrium state for the Boltzmann
operator. In this way a dissipative Euler system is derived.

\end{document}